
\input amstex
\NoBlackBoxes
\documentstyle{amsppt}
\hoffset=1cm
\pagewidth{138mm}
\voffset=1cm
\pageheight{220mm}
\topmatter

\def\wt{\widetilde}

\title
The determination of integral closures and geometric applications
\endtitle
\author
Sheng-Li Tan and De-Qi Zhang
\endauthor

\abstract We express explicitly the integral closures of some ring
extensions; this is done for all Bring-Jerrard extensions of any
degree as well as for all general extensions of degree $\le 5$; so
far such an explicit expression is known only for degree $\le 3$
extensions. As a geometric application we present explicitly the
structure sheaf of every Bring-Jerrard covering space in terms of
coefficients of the equation defining the covering; in particular,
we show that a degree-$3$ morphism $\pi : Y \rightarrow X$ is
quasi-etale if and only if $c_1(\pi_*{\Cal O}_Y)$ is trivial
(details in Theorem 5.3). We also try to get a geometric
Galoisness criterion for an arbitrary degree-$n$ finite morphism;
this is successfully done when $n = 3$ and less satifactorily
done when $n = 5$.
\endabstract

\endtopmatter
\address
S.-L. Tan: Department of Mathematics, East China Normal University, Shanghai 200062,
P. R. of China. {\hskip1cm \eightrm Email: \ sltan\@math.ecnu.edu.cn}
\endaddress
\address {\rm Current Address:}
Department of Mathematics, National University of Singapore,
2 Science Drive 2, Singapore 117543, Singapore.
{\hskip1cm \eightrm Email: \ mattansl\@nus.edu.sg}
\endaddress
\address D.-Q. Zhang:
Department of Mathematics, National University of Singapore,
2 Science Drive 2, Singapore 117543, Singapore. {\hskip1cm \eightrm Email: \
matzdq\@nus.edu.sg}
\endaddress

\document

\head
Introduction
\endhead

The computation of the integral closure (or normalization)
of a finite extension is a fundamental problem in number theory, commutative algebra
and algebraic geometry. David Hilbert's determination of the ring of
algebraic invariants is exactly the computation of normalization (see [St]).
Readers will also notice that Noether's normalization
theorem appears as the basic pre-requisite of Algebra at page 3 of
David Mumford's "The red book of varieties and schemes".

\par
The computation of integral closure for cyclic extensions is
well known (Lemma 1.3). For a cubic extension, this computation
was pioneered by A. A. Albert in the 1930's [Al1, Al2], and continued and completed
half century later in [Sp], [ShS] over {\bf Z} and in [Ta1] over a
Noetherian unique factorization domain (UFD for short).
They extended a result of Richard Dedekind [De] published in 1899,
who has given an integral basis for a pure cubic field.
However, as far as the authors know, if the extension is non-cyclic
with degree higher than 3, it seems that no general formula has been found.

\par
For a general affine domain, important pioneering works have been done
by W. V. Vasconcelos in [Va1] who established a very
effective algorithm to compute the integral closure (see also [BV] and [dJ]).
For the recent development, the analysis of the algorithm and its
complexity and more background, we refer readers to
Vasconcelos's very comprehensive book [Va2] as well [Va3] and [SUV].
The results in the book [Va2] are very important in simplifying our argument
and make our description of integral closure possible.

\par
A finite extension $R[\alpha]$ of a \text{\rm UFD} $R$ is given by a root $\alpha$
of an irreducible polynomial in $R[z]$:
$$
f=z^n+a_{n-1}z^{n-1}+\cdots+a_1z+a_0, \hskip0.3cm a_i\in R.
$$
By a linear Tschirnhaus transformation $z\mapsto z-a_{n-1}/n$, one
can assume that the coefficient of $z^{n-1}$ vanishes. Jerrard
proved that by a Tschirnhaus transformation involving square and
cube roots, the second, third and fourth terms (after the leading
term) of a general polynomial $f$ can be removed. This result
generalized Bring's result for quintic polynomials. (See [Dh,
pp.195--200]).

\par
To be precise, a Tschirnhaus transformation is the substitution $w
= \alpha_0 + \alpha_1 z + \cdots + \alpha_m z^m$ for some $m \le
n-1$ and some $\alpha_i$ to be determined by solving square and
cube equations with solutions in an over ring ${\hat R}$ of $R$, 
so that if $f(z) = 0$ then $g(w) = 0$ for some
polynomial $g$ in ${\hat R}[w]$ of the form
$$g = w^n + b_{n-4} w^{n-4} + \cdots + b_1 w + b_0.$$
In particular, when $n \le 5$, the determination of roots of
a general degree $n$ polynomial
can be reduced to that of a polynomial of the following form
$$
z^n+sz+t
$$
after making one base change of degree 2 if $n=4$ or
some base changes of degree $\le 3$ if $n=5$. Such a
polynomial is called a {\it Bring-Jerrard polynomial} (or simply
a B-J {\it polynomial}).
The corresponding extension (with the polynomial irreducible)
is called a {\it Bring-Jerrard extension} (or simply a B-J {\it extension}).

\par
The purpose of this paper is to express explicitly the
integral closure $\widetilde A$ of a Bring-Jerrard extension
$A=R[\alpha]$ given by the root $\alpha$ of a degree $n$ Bring-Jerrard polynomial
over $R$ (Theorem 2.1). As an application, we calculate explicitly
the integral closure of a general degree $4$ or $5$ extension by
reducing it to a type B-J extension, using Tschirnhaus transformation
(Theorem 3.1 and Theorem 4.5).

\par
We compute also the ramification divisor of a B-J extension (Theorem 1.4).
Late, we apply the computation of integral closure to algebraic geometry;
especially, we determine explicitly the structure sheaf of the covering space
in terms of coefficients of the equation defining the covering (Theorem 5.1).
As a further application, we prove a geometric criterion for a
degree-$3$ finite morphism to be Galois (Theorem 5.5);
the general degree$-n$ case is more complicated, and we prove
a partial result for the case of degree-$5$ (Proposition 5.7).

\par
The idea of our computation is as follows. The normalization
ring $\widetilde A$ must be a reflexive $R$-module for any
$f$ [Ha2]. On the other hand, any reflexive $R$-module
$M$ of rank $n$ is a syzygy module
$$
0\to M \to R^{m+n} \ {\overset \varphi_M \to \longrightarrow } \ R^{m}
$$
for some $m$, where $\varphi_M$ is an $(m+n)\times m$ matrix over $R$.
In order to find the integral closure ${\wt A}$ defined by $f$, we need to find
the matrix $\varphi_f=\varphi_{\widetilde A}$ from $f$.
This reduction has the advantage that
the syzygy module $M = \widetilde A$ is always reflexive, and so satisfies
automatically Serre's condition $S_2$ (see [Va1] or [Va2]). Hence we only
need to compute the co-dimension one normalization. Serre's condition $R_1$
(co-dimension one nonsingular) is invariant under localization $R_{\frak p}$
at all height-one prime ideals $\frak p$ and under completion
$\widehat R_{\frak p}$. Hence we only need to compute the normalization of
an extension over the one-dimensional ring $\widehat R_{\frak p}$.
By Cohen's Structure Theorem for regular rings
[Ha1, p.34], $\widehat R_{\frak p}\cong k[[x]]$
(assuming $R$ contains a field).
Thus we can assume that $R$ is the ring of formal power series
over the residue field $k$. The polynomial $f$ defines a
local curve in $\bold A_k^2$. Now the normalization is just the
resolution of plane curve singularities, which can be realized by
the embedded resolution [Ha1, p.391].
The next step is to globalize the local computation to $R$.

\par
In fact, the syzygy presentation of the integral closure is the
simplest one, because if $R$ is a general \text{\rm UFD} but not a
PID, then there is no canonical method to solve the syzygy
equations; namely, it is hopeless to give explicitly the
generators of $\widetilde A$ when $n\geq 3$. Indeed, Akizuki [Re,
\S 9.5] shows that there is a 1-dimensional Noetherian local ring
$R$ such that the integral closure ${\wt R}$ of $R$ in its field
of fraction is not finitely generated as an $R$-module; in other
word, the number of generators of ${\wt R}$ as an $R$-module is
not finite, or the syzygy equations in the description of ${\wt
R}$ can not be solved.

\par
The result in this note is also related to the works of Catanese [Ca], Miranda [Mi],
Casnati - Ekedahl [CE] and Hahn - Miranda [HM].
Compared with their results, our approach emphasizes more
on the very close relation between the structure sheaf of
the variety upstairs and the coefficients of the equation (downstairs)
defining the finite morphism.

\head
Assumption
\endhead
For integral domain $R$ we assume that $\text{\rm Char} \, R$ is coprime to
both $n$ and $n-1$, where $n$ is the degree of the extension
$R \subset R[\alpha]$.

\head
Acknowledgment
\endhead
The authors would like to thank Professor F. Catanese for the
interest in the work and the advice on degree-4 covers, and to
thank the referee for useful suggestions. This work was started
during the first named author's visit to National University of
Singapore between October 2001 - March 2002. The first author is
supported by 973 Science Foundation, Doctoral Program Foundation
of EMC and the Foundation of Shanghai for Priority Academic
Discipline. The second author is supported by an Academic Research
Fund of National University of Singapore.

\head
1. Discriminant and ramification of a Bring-Jerrard extension
\endhead

Let $n \ge 3$ be an integer. In this section, we shall calculate explicitly
the discriminant and ramification of the integral closure of
a Bring-Jerrard extension $R \subset R[\alpha]$ of a Noetherian \text{\rm UFD} $R$
containing a field, which is defined by a root $\alpha$ of $f(z) = z^n+sz+t$
with non-zero elements $s,t$ in $R$.

\par
Following [Ta1], we shall decompose $s$ and $t$ as well as
the discriminant of $f$ as the products
of elements defining reduced divisors in ${\Cal Spec} \, R$.
For a prime element $p$ in $R$, we let $s_p=\nu_p(s)$ be the corresponding
valuation of $s$. Now we define:
$$
\gather
\varepsilon_p:= n t_p - (n-1) s_p, \hskip0.5cm
\lambda_p:= \min\left\{
\left[ \frac{s_p}{n-1} \right],
\left[ \frac{t_p}{n} \right]
\right\}, \\
a_i := \underset
{\varepsilon_p \equiv i \ (n)} \to
{\prod_{\varepsilon_p>0}} p ,
\hskip0.5cm
b_j :=
\underset
{\varepsilon_p \equiv j \ {(n-1)}} \to
{\prod_{\varepsilon_p< 0}} p,
\endgather
$$
where $1 \le i \le n-1$ and $1 \le j \le n-2$.
If $\lambda_p > 0$ then $f/p^n = (z/p)^n + (s/p^{n-1}) (z/p) + (t/p^n)$
is in $R[z]$ and we may replace $f$ by $f/p^n$. Thus we may and will assume
that $\lambda_p = 0$ for all prime element $p$, i.e., the data $(s, t)$ is minimal
(see {\bf (1.1)} below). We can check easily the following decomposition:
$$
s = a_0 \prod_{i=1}^{n-1} a_i^i
\prod_{j=1}^{n-2} b_j^{n-1-j}, \hskip0.3cm
t = b_0 \prod_{i=1}^{n-1} a_i^i
\prod_{j=1}^{n-2} b_j^{n-j}.
$$

\par
The usual discriminant $\delta$ of $f(z)$ is
given as follows (or the negative of it):
$$\align
\delta &= (n-1)^{n-1} s^n - (-n)^n t^{n-1} \cr
&= \left(\prod_{i=1}^{n-1} a_i^{i}\right)^{n-1}
\left(\prod_{j=1}^{n-2} b_j^{n-1-j}\right)^n c,
\endalign
$$
where
$$
c := a + b
$$
and
$$ a := (n-1)^{n-1} a_0^n \prod_{i=1}^{n-1} a_i^i, \hskip0.5cm
b := -(-n)^n b_0^{n-1} \prod_{j=1}^{n-2} b_j^j.
$$

\par
Let $c_1 := \prod_{c_p=\text{\rm odd}} \, p$.
Then we can write $c = c_0^2 c_1$. From the definition of
$a_i, b_i$ and the relation $c = a+b$, we see easily:
\roster
\item"$\bullet$" Since $s\neq 0$, these $a$, $b$ and $c$ are relatively coprime;
if $i \ge 1$, then $a_i$, $b_i$ and $c_1$ are square-free,
i.e., they define reduced divisors of ${\Cal Spec} \, R$;
\item"$\bullet$" $a_0$, $b_0$ and $c_0$ are not necessarily square-free; and
$$\delta =
\left({\prod_{i=1}^{n-1}  {a_i}^i}\right)^{n-1}
\left({\prod_{j=1}^{n-2} {b_j}^{n-1-j}}\right)^n
c_1  c_0^2.$$
\endroster

\par
\definition{Definition 1.1} (1) $z^n+sz+t$ or the data $(s, t)$
is called {\it minimal} if there is no
prime divisor $p$ such that $p^{n-1}\,|\,s$ and $p^n\,|\, t$, i.e., $\lambda_p=0$
for all prime divisors $p$.

\par
(2) $z^n+sz+t$ is said to be equivalent to $z^n+s'z+t'$ if there is a
regular section $e$ on $X$ without zero point on $X$ such that
$s'=e^{n-1}s$ and $t'=e^nt$.

\par
(3) The triplets $(a,b,c)$ and $(a',b',c')$ of coprime sections of a
line bundle with $a+b=c$ and $a'+b'=c'$
are said to be equivalent if there is a regular section
$e$ without zero point on $X$ such that $a'=ea$, $b'=eb$ and $c'=ec$.
\enddefinition

\par
If the data $(s, t)$ is not minimal, i.e., if $\lambda=\prod_{p}p^{\lambda_p}$
is not a unit, then we can check easily that the defining polynomial
of $\alpha/\lambda$ is a minimal B-J polynomial of degree $n$,
and the normalization ring of $R[\alpha/\lambda]$ is equal to that
of $R[\alpha]$. So we can always assume that the data $(s, t)$ is minimal.
Clearly, we have:

\proclaim{Proposition 1.2}
Let $R$ be a Noetherian \text{\rm UFD}.
Then up to equivalences {\bf (1.1)},
there is a one to one correspondence between the following two sets:
$$
\left\{\text{ \rm Minimal } z^n+sz+t \text{ \rm with } s\neq 0\,\right\}
\longleftrightarrow
\left\{\text{ \rm Coprime } (a,b,c) \text{ \rm with } a+b=c \,
\right\}.
$$
\endproclaim

\par
The following result is well-known and stated for later use (see
[EV, pp.18--35]).

\proclaim{Lemma 1.3}
Let $R$ be a Noetherian \text{\rm UFD}, and let $A=R[\alpha]$ be a
cyclic extension of $R$ which is defined by a root $\alpha$ of
an irreducible polynomial $z^n + u\ell_1 \ell_2^2 \dots \ell_n^n$ in $R[z]$,
where $\ell_1,\cdots,\ell_{n-1}$ in $R$
are square-free and $u$ in $R$ is a unit. Denote by $\widetilde A$
the integral closure of $A$. Then we have:
\roster
\item As an $R$-module, $\widetilde A$ is generated by:
$$ {\alpha^k\over
\prod_{j=1}^{n} \ell_{j}^{[jk/n]}}, \hskip0.3cm k=0,1,\cdots,{n-1}.$$
\item Denote by $\frak p$ a minimal ideal of $R$ containing $\ell_k$,
and by $\widehat R_{\frak p}$ the completion of the local ring $R_{\frak p}$.
Then the completion $\widetilde A \otimes_{R_{\frak p}}\widehat R_{\frak p}$
of $\widetilde A$ has exactly $\gcd(n,k)$ maximal ideals (all of height $1$)
whose generators ${\bar x}_{ki}$ can be chosen so that we have the following factorization:
$$
\ell_k = \left(\prod_{i=1}^{\gcd(n,k)}\bar x_{ki}\right)^{n\over\gcd(n,k)}.
$$
\endroster
\endproclaim

\par
Now we can prove the main result of this section.

\proclaim{Theorem 1.4}
Let $n \ge 3$ be an integer.
Let $R$ be a Noetherian \text{\rm UFD} containing a field with
$\text{\rm Char} \, R$ coprime to $n$ and $n-1$, and let $A=R[\alpha]$ be a
$B-J$ extension of $R$ given by a root $\alpha$ of an irreducible polynomial
$f(z) = z^n+sz+t$ in $R[z]$ with $s \ne 0$
and the data $(s, t)$ minimal (see ${\bold (1.1)}$).
Denote by $\widetilde A$ the integral closure of $A$. Then we have:
\roster
\item The discriminant of the extension $\widetilde A$ over $R$
(i.e., the defining equation of the branch locus of the finite cover
${\Cal Spec} \, {\wt A} \rightarrow {\Cal Spec} \, R$) is
$$
D_{\widetilde A/R}=c_1\prod_{k=1}^{n-1}a_k^{n-\gcd(n,k)}
\prod_{k=1}^{n-2}b_k^{n-1-\gcd(n-1,k)}
$$
\item
Let $\frak p$ be a prime ideal of $R$ generated by a factor $x$ of $a_k$ or $b_k$ or $c_k$.
Then we have the following factorizations into primes in the completion
$\widetilde A\otimes _{R_{\frak p}}\widehat R_{\frak p}$ of $\widetilde A$,
$$x=\cases
\left(\prod_{i=1}^{\gcd(n,k)}{\bar x_i}\right)^{\ n \over\gcd(n, k)},
&\text{ if } \ x\,| \, a_k,  \\
\left(\prod_{i=1}^{\gcd(n-1,k)}{\bar x_i}\right)^{\ n-1\over \gcd(n-1, k)} \, {\bar x_1}',
&\text{ if } \ x\,| \, b_k , \\
{\bar x_1}^2 \ \prod_{i=1}^{n-2}{\bar x_i}', &\text{ if } \ x\,| \, c_1. \\
\endcases
$$
In particular, we have the following factorizations of $a_k$, $b_k$ and $c_1$ in
$\widetilde A$, modulo the multiplication by a unit of $\wt A$,
where $\bar a_k$, $\bar b_k$ and $\bar c_1$ ($k \ge 1$) are square free
$$
a_k ={\bar a_k}^{\ n/\gcd(n, k)},  \hskip0.3cm
b_k ={\bar b_k}^{\ (n-1)/\gcd(n-1, k)} \ {\bar b_k}',
\hskip0.3cm
c_1 = {\bar c_1}^2 \ {\bar c_1}'.
$$
\item The defining equation of the ramification
divisor of the finite cover ${\Cal Spec} \, {\wt A} \rightarrow {\Cal Spec} \, R$ is
$$
\Cal R_{\widetilde A/R} = \bar c_1
\prod_{k=1}^{n-1}
{\bar a_k}^{\ {n\over \gcd(n,k)} - 1}
\prod_{k=1}^{n-2}
{\bar b_k}^{\ {n-1\over \gcd(n-1,k)} - 1}.
$$
\endroster
\endproclaim

\demo{Proof} Note that the discriminant
$D_{\widetilde A/R}$ is a factor of
$$
\delta=D_{A/R}
= (n-1)^{n-1} s^n - (-n)^n t^{n-1}
=\left({\prod_{i=1}^{n-1}  {a_i}^i}\right)^{n-1}
\left({\prod_{j=1}^{n-2} {b_j}^{n-1-j}}\right)^n
c_1  c_0^2.$$
So we need only to find the ramification indices over the prime ideal
$\frak p=(x)$ generated by a factor $x$ of $a_k$ ($k \ge 1$), $b_k$ ($k \ge 1$)
or $c_k$.

\par
Localizing $R$ at $\frak p$,
we may assume that $R=R_{\frak p}$, which is a DVR, so that
$x$ is a parameter. Denote by $\frak m$ a (height 1) minimal ideal
of $\widetilde A$ over $\frak p$. Since we consider only the ramification
at $\frak m$ over $\frak p$, we may reduce to the completion $\widehat R$
with respect to $\frak p$. So we may assume that $R$ is complete.
By Cohen's Structure Theorem [Ha1, p.34],
a complete regular local ring of dimension 1 containing some field is the
ring of formal power series over the residue field $k$, i.e.,
$R\cong k[[x]]$.
Hence $A$ is a local curve over $k$ defined by $z^n+s(x)z+t(x)=0$
in $\bold A_k^2$. $\widetilde A$ is the normalization of this curve.
The ramification index of $\widetilde A$ at $\frak m$ over $\frak p$
is equal to the corresponding ramification index as a local $n$-cover
over ${\Cal Spec} \,(k[[x]])$.
Since the normalization can be realized by the embedded resolution
of a plane curve singular point, in what follows, we shall
compute the ramification index by using Lemma 1.3 and the embedded
resolution.

\par
Note that $\widetilde A$ is a finite module over $R$, so $\widetilde A$
is also complete.

\par
We consider first the case $x\,|\,a_k$. We rewrite the polynomial:
$$
f(z) = z^n + a' x^m z + b' x^k = z^n + (b' + a' x^{m-k} z)x^k=z^n+ux^k,
$$
where $m \ge k$, $u=b'+a'x^{m-k}z$, and $a'$ and $b'$ are units in $R$.
Since $f(\alpha) = 0$ and $x$ is in $\frak m$, we see that $z = \alpha$ is in $\frak m$.
Thus $u$ is a unit of $\widetilde A$. So we are reduced to Lemma 1.3.
This proves Theorem 1.4 for factors of $a_k$.

\par
Next we consider the case $x\,|\, b_{n-1-k}$.
So
$$
f = z^n + a' x^k z + b' x^m
$$
where $m \ge k+1$, and $a'$ and $b'$ are units in $R$.
Now we blow up $\bold A_k^2$ at $(0,0)$. Then the strict transform
of the curve $C: f(z) = 0$, is locally a union of an irreducible component
isomorphic to the original curve $C$ and a curve (of degree $n-1$ over $C$) given below,
where $z$ and $x' = x/z$ are new coordinates and the divisor $x' = 0$ is also
the strict transform of the divisor $x = 0$ on $C$
$$
z^{n-k-1}+{x'}^k v = 0, \,\,\,\, v := a'+b'{x'}^{m-k}z^{m-k-1}
$$

\par
Since this $v$ is a unit, this case is
reduced to Lemma 1.3 in applying which we note also that the total transform in
${\Cal Spec} \, {\wt A}$ of the divisor $x = 0$ in ${\Cal Spec} \, A$ is defined
by $zx' = 0$. This proves Theorem 1.4 for factors of $b_{n-1-k}$.

\par
Finally, we consider the case $x\,|\,c_k$ for $k = 0$ or $1$.
In this case, $s$ and $t$ are units of $R$.
Set $w := z + {n t\over (n-1)s}$. It is easy to know
that if there is a ramification over $x=0$, then
we must have $w=0$. We rewrite the polynomial
$$
g(w) := f = w^n + \sum_{i=0}^{n-1} e_i w^i,
$$
where
$$e_i=\cases
-\dfrac{t \delta}{((n-1)s)^n}, & i=0,  \cr
\dfrac{\delta}{((n-1)s)^{n-1}}, & i=1, \cr C_n^i \left(\dfrac{-n
t}{(n-1)s}\right)^{n-i}, &i \ge 2.
\endcases
$$
Hence $e_i$ is a unit when $i\geq 2$.
Clearly there is an $m \ge 1$ such that the following are all units in $R$
$$\delta/x^m, \,\, c_0^2c_1/x^m, \,\, e_i/x^m \,\, (i = 0, 1).$$
We rewrite $g(w)$ as follows, where $a', b'$ are units of $R$ and
$u_1$ and $u_2$ are units of $\widetilde A$.
$$\align
g(w) &= w^n + \sum_{i=2}^{n-1} c_i w^i + a' x^m w + b' x^m\cr
     &= w^2u_1+x^mu_2.
\endalign
$$
We see easily that if $x \, | \, c_1$
then $m$ is odd and the ramification
index is 2; if $x$ does not divide $c_1$ (and hence $x$ divides $c_0$)
then the normalization has no ramification. This proves Theorem 1.4
for factors of $c_k$. The proof of Theorem 1.4 is completed.
\enddemo

\head
2.  Integral closure of a B-J extension
\endhead

In this section, we will calculate explicitly the integral closure
of a B-J extension.

\par
Let $n \ge 3$ be an integer.
Let $R$ be a Noetherian \text{\rm UFD} with $\text{\rm Char} \, R$ coprime to $n$ and $n-1$, and let $A=R[\alpha]$ be a
B-J extension of $R$ defined by a root $\alpha$ of an irreducible
polynomial $f(z) = z^n+sz+t$ in $R[z]$ with $s \neq 0$.
Suppose that the data $(s, t$) is minimal,
i.e., $\lambda_p = 0$ for all prime divisors $p$; see {\bf (1.1)}.
Denote by $\widetilde A$ the integral closure of $A$, and by $\widetilde A_0$
the trace-free $R$-submodule of ${\wt A}$. Obviously the trace map
$\text{Tr}: \ \widetilde A\to R$ splits because $\text{\rm Char} \,  R$
is coprime with $n$ (and also $n-1$) by the assumption in the Introduction.
Thus
$$
\widetilde A=R\oplus \widetilde A_0.
$$
For $1 \le i \le n-1$, we set
$$
h_i := \prod_{k=1}^{n-1} a_k^{[ki/n]}
\prod_{k=1}^{n-2} b_k^{[(n-1-k)i/(n-1)]}.
$$
If we denote by $\sigma_i$ the coefficient of $z^{n-i}$ in $f$,
then we can compute $s_i=\text{Tr}(\alpha^i)$ by Newton's identities:
$$\align
&s_1+\sigma_1=0,\cr
&s_2+s_1\sigma_1+2\sigma_2=0,\cr
&\hskip1cm \vdots\cr
&s_k+s_{k-1}\sigma_1+\cdots+s_1\sigma_{k-1}+k\sigma_k, \hskip0.3cm k\leq n.
\endalign
$$
In our case, $\sigma_i=0$ for $i<n-1$,
$\sigma_{n-1}=s$ and $\sigma_n=t$.
Thus we get
$$
\text{Tr}(\alpha^i)=0,  \text{ for } i=1,\cdots, n-2, \
\text{ and } \text{Tr}(\alpha^{n-1})=-{(n-1)s}.
$$
So we can construct $n-1$ trace free elements:
$$
\beta_i:=\cases \dfrac{\alpha^i}{ h_i}, &\text{ if } \ i=1,\cdots,n-2, \cr
\dfrac{\alpha^{n-1}+{n-1\over n}s}{ h_{n-1}}, &\text{ if } \ i=n-1.
\endcases
$$
For $i=1,\cdots,n-2$, we define
$$\align
f_i:&=(n-1)a_0\prod_{k=1}^{n-1}a_k^{\left[\frac{(i+1)k}{n}\right]-\left[\frac{ik}{n}\right]} ,\cr
g_i:&=nb_0\prod_{k=1}^{n-2}b_k^{1+\left[\frac{(n-1-k)i}{n-1}\right]-\left[\frac{(n-1-k)(i+1)}{n-1}\right]}.
\endalign
$$

\par
We now state the main result of this section.

\proclaim{Theorem 2.1} Let $n \ge 3$ be an integer. Let $R$ be a
Noetherian \text{\rm UFD} containing a field with $\text{\rm Char}
\, R$ coprime to $n$ and $n-1$ and let $A = R[\alpha]$ be a $B-J$
extension given by a root $\alpha$ of an irreducible polynomial
$f(z) = z^n + s z + t$ in $R[z]$ with $s \ne 0$ and the data $(s,
t)$ minimal. Then the integral closure ${\wt A}$ of $A$ in the
fraction field of $A$ satisfies ${\wt A} = R \oplus {\wt A}_0$, as
$R$-modules, where the trace $($over $R)$ free $R$-submodule ${\wt
A}_0$ of ${\wt A}$ is given as follows:
$$
\widetilde A_0=\left\{\left. \dfrac{v_1\beta_1+v_2\beta_2+\cdots+v_{n-1}\beta_{n-1}}{c_0}
\,\right|\, v_j\in R, \ c_0\,|\,f_kv_k+g_kv_{k+1}, \ 1\leq k\leq n-2 \,\right\}.
$$
\endproclaim

\proclaim{Corollary 2.2} With the assumptions in ${\bold (2.1)}$,
$\widetilde A_0$ is the following syzygy module
$$
0\to \widetilde A_0 \to R^{2n-3} \ {\overset M\to\longrightarrow} \ R^{n-2}.
$$
Here $M$ sends $(v_1,\cdots,v_{2n-3})$ to
$(f_1v_1+g_1v_{2}+c_0v_{n},
f_2v_2+g_2v_{3}+c_0v_{n+1},\cdots,
f_{n-2}v_{n-2}+g_{n-2}v_{n-1}+c_0v_{2n-3})$.
\endproclaim

\demo{Proof of Theorem $2.1$} Note first that the $\widetilde A$ and $\widetilde A_0$
are reflexive [Ha2]. On the other hand, the right hand side of the displayed
equality in the theorem is a syzygy module, so it is also reflexive. In order
to get the equality, we only need to prove that the syzygy is the
co-dimension one normalization of $A$. Because being reflexive implies the $S_2$
condition in Serre's Criterion [Ha1, p.185], or two reflexive modules over $R$
are isomorphic if and only if they are co-dimension 1 isomorphic.
Clearly, the right hand side of the displayed equality in the theorem, added with the
summand $R$, contains $R[\alpha]$.

\par
Now we use the same technique as in the proof of Theorem 1.4. Namely,
we reduce the proof to the case when $R=k[[x]]$. Hence
$A=k[[x,z]]/(z^n+s(x)z+t(x))$. Then we only need to prove that
the normalization ring of the local curve singularity defined by
$f=z^n+s(x)z+t(x)=0$ can be generated by the syzygies in the theorem.
Clearly, the theorem is true outside the set $\delta = 0$.
So we may assume that $x\,|\,\delta$ and only need to check the
equality in the theorem at the smooth points of Supp$(\delta = 0)$.

\par
If $x\,|\,a_k$, then as \S1, we have $f=z^n+ux^k$, where $u$ is a unit.
Using Lemma 1.3, we can see easily that the normalization ring
is generated by $\beta_i$. Hence it is generated by the syzygies.

\par
If $x\,|\, b_{n-1-k}$, then
$$
f = z^n + a' x^k z + b' x^m
$$
where $m \ge k+1$, and $a'$ and $b'$ are units in $R$.

\par
We shall prove first that these $\beta_i$ are integral over $R$.
Since the localization $\widetilde A_{\frak m}$ of $\widetilde A$
at a height-1 prime ideal $\frak m$ over $(x)$ is a DVR, we can define a valuation
$\nu=\nu_{\frak m}$. Thus for any two elements $g$ and $h$ in $\widetilde A$,
the element $r=g/h$ is in $\widetilde A_{\frak m}$ if and only if
$\nu(r)=\nu(g)-\nu(h)\geq 0$.
Now we claim that $\nu(\beta_i)\geq0$ for all $i$.

\par
Suppose the contrary that $\nu(\beta_i)<0$ for some $i$.
If $i<n-1$, then we have
$$
i\nu(\alpha)<\nu(h_i)=\left[{ki\over n-1}\right]\nu(x)\leq {ki\over n-1}\nu(x),
$$
so $\nu(\alpha^{n-1}) = (n-1)\nu(\alpha) < k \nu(x)$ and $\nu(\alpha^n) < \nu(x^k \alpha)$.
This and the equation
$\alpha^n + a' x^k \alpha + b' x^m=0$ imply that $\nu(\alpha^n)$ must be equal
to $\nu(x^m)$. Hence
$$
{m\over n}\nu(x) = \nu(\alpha) <{k\over n-1}\nu(x).
$$
This implies that ${n\over n-1}>{m\over k}\geq {k+1\over k}$. Hence we have
$k>n-1$, a contradiction. Thus $i=n-1$. Then we have
$$
\nu\left(\alpha^{n-1}+{n-1 \over n}s\right)=
\nu\left(\alpha^{n-1}+{n-1 \over n}a'x^k\right)<\nu(h_{n-1})=\nu(x^k)=\nu(s),
$$
so $\nu(\alpha^{n-1})<\nu(s)=k\nu(x)$, a contradiction as above.
Hence $\beta_i$ must be integral over $R$ for all $i$.

\par
Now we need to prove that these $\beta_i$ generate the normalization ring. It is
enough to prove that the discriminant $d=\det(\text{Tr}(\beta_i\beta_j))$
of $\beta_0=1,\beta_1,\cdots, \beta_{n-1}$
is equal to the discriminant $D_{\widetilde A/R}$ (see \S1).
From the definition of $\beta_i$, we get the following,
where we set $h_0:=1$ and replace $\beta_{n-1}$ by $\alpha^{n-1}/h_{n-1}$
$$
d=\det\left(\text{Tr}\left({\alpha^i\alpha^j\over h_ih_j}\right)\right)
={\det(\text{Tr}(\alpha^i\alpha^j))\over\prod_{i=1}^{n-1}h_i^2}={\delta\over \prod_{i=1}^{n-1}h_i^2}.
$$
The identity $d=D_{\widetilde A/R}$ (up to a unit) is equivalent to
$$
2\sum_{i=1}^{n}\left[{ki\over n} \right]=(n+1)k-n+\gcd(n,k),
$$
or
$$
{2\over n}\sum_{i=1}^{n}\varepsilon_i-n+\gcd(n,k)=0,
$$
where $0 < k < n$ and $\varepsilon_i=ki-n[ki/n]$ is a non negative integer
less than $n$. Note that if $\gcd(n,k)=1$, then
$\varepsilon_i$ attains all of the numbers between 0 and $n-1$,
thus it is easy to obtain the above identity.
The proof for the general case can be reduced to this case.
So we have proved that these $\beta_i$ are the generators of the
normalization ring.

\par
Finally we consider the difficult case when $x\,|\,c_1c_0^2$.
As in \S 1, we rewrite $f(z)$ as follows, where $e$ is a unit near $w = 0$
$$
f(z)=w^2e+we_1+e_0.
$$
Since $c_0^2c_1 \, | \, e_i$ ($i = 0, 1$) (\S 1), we see from the above equation that
$we/c_0$ is integral over $A=R[\alpha]$, so it is also integral
over $R$.

\par
We shall find $u_i$ in $R$ such that
$$
\beta':={we\over c_0}-
{1\over n}{\text{Tr}\left({we\over c_0}\right)}=
{u_1\beta_1+\cdots+u_{n-1}\beta_{n-1}\over c_0}.
$$
If we let $q=-nt/(n-1)s$, then
$$
we=\sum_{i=1}^{n-1}q^{n-1-i}z^i-(n-1)q^{n-1}.
$$
Thus it is easy to see that
$$
u_i=q^{n-1-i}h_i.
$$

\par
Now we claim that $\beta=\sum_{i=1}^{n-1}v_i\beta_i/c_0$
is integral over $R$ if $v_i$ satisfies the conditions
in the expression of ${\wt A}_0$ in the theorem.
In fact, one can check that these conditions
are equivalent to the following
$$
c_0 \ | \ u_{i+1}v_i-u_iv_{i+1}, \hskip0.3cm i=1,\cdots,n-2.
$$
In particular, $\beta'$ satisfies the conditions.
We consider the element
$\beta''=u_1\beta-v_1\beta'=(\sum_{i=1}^{n-1}v_i''\beta_i)/c_0$,
$v_1''=0$. By using the induction on $i$ and the above conditions,
we see that $v_i''$ is divided by $c_0$. Hence $\beta''$ is
integral over $R$. Since $u_1$ is a unit, $\beta$
is also integral over $R$.

\par
Next we need to prove that the syzygies generate the normalization ring ${\wt A}$
as an $R$-module. In fact, we only need to prove that
$\beta^*=we/c_0$ generates the normalization ring
near $w = 0$ (see \S 1). We shall show that $R[\alpha, \beta^*]
= R[w, \beta^*] = R + R \beta^*$ is normal near $x = w = 0$.
As in \S 1, let $m \ge 1$ such that $c_0^2c_1/x^m, e_i/x^m$ ($i = 0, 1$) are units in $R$
near $x = 0$. Here
$m$ is odd if and only if $x \, | \, c_1$.
If $m = 1$ then the equation $f(z) = w^2e + we_1 + e_0 = 0$ has a linear term
(near $x = w = 0$) and hence $R[\alpha] = R[z]/(f(z))$ is smooth near $w = 0$.
We may assume $m \ge 2$ and hence $m \, | \, c_0$ for $c_1$ is reduced.
Note that $R[w, \beta^*]$ is the quotient of $R[w, y] = k[[x]][w, y]$
modulo the equations below, where $e$ is a unit in $R$ near $w = 0$
$$F_1 := w - yc_0/e = 0, \,\,\,\, F_2 := y^2+{e_1\over c_0}y+{e_0\over c_0^2}e = 0.$$

\par
If $x \, | \, c_1$ then $F_2 = 0$ implies that $y = 0$ (when $x = 0$)
and $F_2$ has a linear term, whence $R[w, \beta^*]$ is smooth near $x = w = 0$.
If $x$ does not divide $c_1$, then the partial derivative $(F_2)_y \ne 0$ holds
near $x = w = 0$ and along the zero locus of $F_2 = 0$,
so is the smoothness of $R[w, \beta^*]$ near $x = w = 0$.
We have completed the proof of Theorem 2.1.
\enddemo

\head
3. Integral closure of a quartic extension
\endhead

In this section, we will calculate explicitly the integral closure of
a quartic extension $R[\alpha]$ of a Noetherian \text{\rm UFD} $R$
with $\text{\rm Char} \, R$ coprime to $2$ and $3$,
which is given by a root $\alpha$ of an irreducible quartic monic polynomial $f(z)$ over $R$.
We may assume, after a shift of coordinate, that
$$
f(z) = z^4+\sigma_2 z^2 -\sigma_3 z+\sigma_4.
$$
\par
We do the factorization:
$$
2\sigma_2^3-8\sigma_2\sigma_4+9\sigma_3^2=d_1d_0^2,
$$
where $d_1$ is square free in $R$, i.e., no square of
a prime element of $R$ divides $d_1$.

\par
Consider the general case where
$d_1$ has no square root in the fraction field of $R[\alpha]$.
So we have ring extensions $R \subset R[\alpha] \subset R[\alpha, y]$
of degrees 4 and 2, where $y^2 = d_1$. Thus the ring extensions
$R \subset R[y] \subset R[\alpha, y]$ are of degrees 2 and 4.
In the new ring ${\hat R} := R[y]$, we can find an element
$$
w = {1\over 2}\sigma_2^2+\left({3\over 2}\sigma_3+{1\over 2}d_0y\right)
\alpha+\sigma_2\alpha^2.
$$
This ${\hat w} = w$ satisfies
$$
\hat w^4+\hat s \hat w+ \hat t=0,
$$
where
$$\align
\hat s = &-{1\over2}\sigma_2^6+4\sigma_2^4\sigma_4-{19\over4}
\sigma_2^3\sigma_3^2-8\sigma_2^2\sigma_4^2+27\sigma_2\sigma_3^2\sigma_4-{27\over2}
\sigma_3^4
\cr
&-{1\over4}yd_0\sigma_3(3\sigma_2^3-28\sigma_2\sigma_4+
18\sigma_3^2),\cr
\hat t =
&{19\over8}\sigma_2^5\sigma_3^2-{5\over4}\sigma_2^6
\sigma_4+\sigma_2^4\sigma_4^2+{3\over16}\sigma_2^8
+{3\over2}\sigma_2^3\sigma_4\sigma_3^2+{81\over2}\sigma_4
\sigma_3^4
\cr
&+{27\over4}\sigma_2^2\sigma_3^4
-36\sigma_4^2
\sigma_3^2\sigma_2+4\sigma_2^2\sigma_4^3\cr
&+yd_0\left({3\over8}
\sigma_2^5\sigma_3+{27\over 2}\sigma_4
\sigma_3^3-6\sigma_4^2\sigma_3\sigma_2+{9\over 4}\sigma_2^2
\sigma_3^3\right).
\endalign
$$

\par
We consider the general case that the above degree-4 polynomial
is irreducible over the  fraction field of ${\hat R}$.
Then we have ring extensions $R \subset {\hat R} \subset {\hat R}[w]$
of degrees 2 and 4. Since ${\hat R}[w] \subseteq R[y, \alpha]$
and both rings are of degree-8 extensions of $R$, they have the same
 fraction field; in particular, they have the same integral closure
${A^*}$ in their common  fraction field $Q(R)[\alpha, y]$.
We denote by $A_0^*$ the ${\hat R}$-submodule of ${A^*}$
consisting of elements of trace zero over ${\hat R}$.
Since $R$ is \text{\rm UFD} and $d_1$ is also square free in $R$, the ring
$\hat R$ is a normal ring, and thus co-dimension-one regular
(see [Ma, \S 9, Example 4, p. 65]).
We can factor any element in $\hat R$ into the product of primes over the
smooth locus of $\hat R$. We do this for $\hat s$ and $\hat t$;
note that $\hat a_i, \hat b_j$ below may have
different expressions at different affine open sets of ${\Cal Spec} \, \hat R$,
and $\text{\rm div}(\hat a_i)$ and $\text{\rm div}(\hat b_i)$ are reduced
divisors for all $i \ge 1$
$$
\hat s=\hat a_0 \, \hat a_1\,\hat a_2^2\,\hat a_3^3\,\hat b_1^2\,\hat b_2,
\hskip0.3cm
\hat t=\hat b_0 \, \hat a_1\,\hat a_2^2\,\hat a_3^3\,\hat b_1^3\,\hat b_2^2.
$$
We can also define $\hat c_k, \hat f_k, \hat g_k, \hat h_k, \beta_k, \cdots$
over the smooth locus of ${\Cal Spec} \, \hat R$ as in \S 2.
By Theorem 2.1, $A_0^*$ is given as follows
(over the smooth locus of ${\Cal Spec} \, \hat R$):
$$
A_0^* =\left\{\left.
\dfrac{\hat v_1\beta_1+\hat v_2
\beta_2+\hat v_3\beta_3}{\hat c_0}
\,\right|\, \hat v_j\in {\hat R}, \ \hat c_0\,|\,\hat f_k\hat v_k+\hat g_k\hat v_{k+1},
\ k=1, \ 2 \,\right\}.
$$
We denote by $\bar \bullet$ the involution of $R[\alpha][y]$ over $R[\alpha]$, i.e.,
$\overline{r_0+r_1y}=r_0-r_1y$ for $r_i$ in $R[\alpha]$.

\par
Let ${\wt A}$ be the integral closure of $R[\alpha]$ in its  fraction field
and let $\widetilde A_0$ be the trace (over $R$) free part.
Note that the following $6$ elements are integral over $R$,
contained in the  fraction field of $R[\alpha]$
(for being the involution $\bar \bullet$-invariant)
and trace (over $R$) free (noting that $\text{\rm Tr} | R[\alpha, y] / R[y]$
is the lifting of $\text{\rm Tr} | R[\alpha] / R$)
$$
\gamma_i={\beta_i+\bar\beta_i\over 2},
\hskip0.3cm
\gamma_{i+3}={\beta_i-\bar\beta_i\over 2}y, \hskip0.3cm
i=1,2,3.
$$
So these $6$ elements are in $\widetilde A_0$. Note that
$$
\beta_i=\gamma_i+\gamma_{i+3}{1\over d_1}y, \hskip0.3cm i=1,2,3.
$$
By the reasoning above, we have:
$$
\widetilde A_0=\left\{\left. {\beta+\bar\beta\over 2}
\, \right | \,\beta\in A_0^* \right\}.
$$
Let
$$\align
&\hat f_k=f_k+f_{k+2}y, \hskip0.3cm \hat g_k=g_k+g_{k+2}y, \hskip0.3cm
k=1,2; \cr
&\hat v_i=v_i+v_{i+3}y, \hskip0.3cm i=1,2,3; \cr
&\hat c_0=c_0+c_0' y.
\endalign
$$
Then the elements of $\widetilde A_0$ can be expressed as
$$
\sum_{i=1}^6v_i\mu_i,
$$
where
$$\align
{\mu _{1}}&:= {c_{0}}\,{\gamma_{1}} - {c_0'}\,{\gamma_{4}},\cr
{\mu _{2}}&:= {c_{0}}\,{\gamma_{2}} - {c_0'}\,{\gamma_{5}},\cr
{\mu _{3}}&:= {c_{0}}\,{\gamma_{3}} - {c_0'}\,{\gamma_{6}},\cr
{\mu _{4}}&:=  - {c_0'}\,{d_{1}}\,{\gamma_{1}} + {c_{0}}\,{\gamma_{4}},\cr
{\mu _{5}}&:=  - {c_0'}\,{d_{1}}\,{\gamma_{2}} + {c_{0}}\,{\gamma_{5}},\cr
{\mu _{6}}&:=  - {c_0'}\,{d_{1}}\,{\gamma_{3}} + {c_{0}}\,{\gamma_{6}}.
\endalign
$$
The 2 syzygies
$$
\align
&\hat f_1 \hat v_1+\hat g_1 \hat v_2+\hat c_0 (v_7+v_9y)=0,\cr
&\hat f_2 \hat v_2+\hat g_2 \hat v_3+\hat c_0 (v_8+v_{10}y)=0
\endalign
$$
induce 4 syzygies $\phi_i=0, \ i=1, \ 2, \ 3, \ 4$, where $\phi_i$
are defined by
$$\align
{\phi _{1}} :&= {f_{1}}\,{v_{1}}
+ {g_{1}}\,{v_{2}}
+ {d_{1}}\,{f_{3}}\,{v_{4}}
+ {d_{1}}\,{g_{3}}\,{v_{5}} + {c_{0}}\,{v_{7}} + {d_{1}}\,{c_0'}\,{v_{9}},\cr
{\phi _{2}} :&= {f_{3}}\,{v_{1}}
+ {g_{3}}\,{v_{2}}
+ {f_{1}}\,{v_{4}}
+ {g_{1}}\,{v_{5}} + {c_0'}\,{v_{7}} + {c_{0}}\,{v_{9}},\cr
{\phi _{3}} :&=
{f_{2}}\,{v_{2}}
+ {g_{2}}\,{v_{3}}
+ {d_{1}}\,{f_{4}}\,{v_{5}}
+ {d_{1}}\,{g_{4}}\,{v_{6}}
+ {c_{0}}\,{v_{8}}
+ {d_{1}}\,{c_0'}\,{v_{10}},\cr
{\phi _{4}} :&=
{f_{4}}\,{v_{2}}
+ {g_{4}}\,{v_{3}}
+ {f_{2}}\,{v_{5}}
+ {g_{2}}\,{v_{6}}
+ {c_0'}\,{v_{8}}
+ {c_{0}}\,{v_{10}}.
\endalign
$$
Next we shall show that there are relations among $\mu_i$.
Indeed, by the generating property of $\beta_1, \beta_2, \beta_3$,
we have
$$({\bar{\beta}}_1, {\bar{\beta}}_2, {\bar{\beta}}_3) =
(\beta_1, \beta_2, \beta_3) M/{\hat c}_0,$$
where $M$ is a $3 \times 3$ matrix with entries in $\hat R$.
Expressing $\beta_i, {\bar \beta}_i$ in terms of $\gamma_j$
we obtain the relation
$$(\gamma_1, \gamma_2, \gamma_3) (\hat c_0 I - M) =
  (\gamma_4, \gamma_5, \gamma_6) (\hat c_0 I + M)/y.$$
On the other hand, the definition of $\mu_i$
implies the following, where $|\hat c_0|^2 = \hat c_0 \bar{\hat c_0}$
$$\align
|\hat c_0|^2(\gamma_1, \gamma_2, \gamma_3) &=
 c_0 (\mu_1, \mu_2, \mu_3) + c_0' (\mu_4, \mu_5, \mu_6), \\
|\hat c_0|^2(\gamma_4, \gamma_5, \gamma_6) &=
 d_1 c_0' (\mu_1, \mu_2, \mu_3) + c_0 (\mu_4, \mu_5, \mu_6).
\endalign$$
So the relation among $\gamma_i$ implies a relation among $\mu_i$:
$$(\mu_1, \mu_2, \mu_3) [\bar{\hat c_0}I - M] =
  (\mu_4, \mu_5, \mu_6) [\bar{\hat c_0}I + M]/y.$$

\par
The main theorem below of this section
follows from the arguments above and the fact that
an element in a normal ring ${\wt A}$ is
determined by its restriction to the open set lying over the smooth locus
of ${\Cal Spec} \, R$ (whose complement has co-dimension at least two).

\proclaim{Theorem 3.1} Let $R$ be a Noetherian \text{\rm UFD}
containing a field with $\text{\rm Char} \, R$ coprime to $2$ and
$3$. For a general degree-$4$ extension $R \subset R[\alpha]$,
where $\alpha$ is a root of an irreducible polynomial $f(z) =
z^4+\sigma_2 z^2 -\sigma_3 z+\sigma_4$ in $R[z]$, the integral
closure $\widetilde A$ of $R[\alpha]$ in the fraction field of
$R[\alpha]$ is given by ${\wt A} = R \oplus {\wt A}_0$, and the
trace $($over $R)$ free $R$-submodule $\widetilde A_0$ of ${\wt
A}$ is given as follows:
$$
\wt A_0=\left\{\,\left.
{\sum_{i=1}^6v_i\mu_i}\,\right|\,
\matrix
v_1,\cdots,v_{10} \,\, \text{\rm in} \,\, R \,\,
\text{form a solution of the}\cr
\text{linear equations} \,\, \phi_k=0,\ k=1,2,3,4\,
\endmatrix
\right\}.
$$
\endproclaim

\remark{Remark 3.2} $(1)$ In Theorem $3.1$, by a general
degree-$4$ extension, we mean that $\tau :=
2\sigma_2^4-8\sigma_2\sigma_4+9\sigma_3^2$ has no square root in
$R$ (so assume it is written as $d_0d_1 = d_0y^2$ with $d_1$
square free in $R$) and that $\hat w^4 + \hat s \hat w + \hat t$
defined above is irreducible over the  fraction field of $R[y]$;
this also implies that $y$ is not in the  fraction field of
$R[\alpha]$, because the extension $R \subset R[y, w]$ is of
degree $8$, whence the extension $R \subset R[y, \alpha]$ is also
of degree $8$ and both extensions share the same  fraction field.

(2) If the $\tau$ above has a square root
in $R$ then the extension $R \subset R[\alpha]$ is very likely to be of type B-J
(at least when ${\hat w}^4 + \hat s \hat w + \hat t$ is irreducible
as a polynomial over $R$). For B-J extension, we refer to Theorem $2.1$.
\endremark

\head 4. Integral closure of a general degree-$n$ extension
\endhead

In this section, we will calculate explicitly the integral closure of a quintic extension $R[\alpha]$ of a Noetherian \text{\rm UFD} $R$ with $\text{\rm Char} \, R$ coprime to $5, 3, 2$, which is given by a root $\alpha$ of an irreducible
quintic monic polynomial $f(z)$ over $R$.

\par
We remark that the case of a general degree-$n$ extension is similar, though the computation
will be more complicated and the simpler polynomial (we may possibly reduce to)
is of the form $z^n + a_{n-4}z^{n-4} + \cdots + a_1 z + a_0$.
We will illustrate by considering the case of degree $n = 5$.

\par
One may assume, after a shift of coordinate, that
$$
f(z) = z^5+\sigma_2 z^3 -\sigma_3 z^2+\sigma_4 z - \sigma_5.
$$

\par \vskip 1pc \noindent
{\bf (4.1) The base change to reduce to type B-J case.}
We shall show that for a suitably general $f(z)$ we can find
$$y := u+v\alpha+w\alpha^2+p\alpha^3+q\alpha^4,$$
where $u, v, w, p, q$ are in a
(relatively not so big) over ring $\hat R$ of $R$,
such that $\hat y = y$ is a zero of the following
type B-J polynomial with coefficients in $\hat R$:
$${\hat y}^5 + {\hat s} {\hat y} + {\hat t}.$$
In other words, after base changes, the extension $R \subset R[\alpha]$ may be
reduced to the extension ${\hat R} \subset {\hat R}[y]$ of B-J type.
We may then apply Theorem 2.1 to get the integral closure of ${\hat R}[y]$
and also that of $R[\alpha]$ in their respective fraction fields.

\par
Let $Z$ denote the matrix representation of the $R$-linear map
$$\alpha : (1, z, z^2, z^3, z^4) \mapsto
(z, z^2, z^3, z^4, z^5 = -(\sigma_2 z^3 -\sigma_3 z^2+\sigma_4 z - \sigma_5)).$$
Then the matrix representation $Y$ of the linear map $y$
is given by $Y = uI+vZ+wZ^2+pZ^3+qZ^4$. Note that $y$ is a zero of
the characteristic polynomial of $Y$:
$$|\lambda I - Y| = \lambda^5 + \sum_{i=0}^4 d_i \lambda^i.$$
Here $d_i$ is a homogeneous polynomial over $R$ of degree $5-i$ in
$u, v, w, p, q$.
We want to find $u, v, w, p, q$ in some over ring ${\hat R}$ of $R$
such that $d_i = 0$ for all $i = 2, 3, 4$. Then we just set ${\hat s} = d_1$
and ${\hat t} = d_0$ and get the desired $y$ satisfying
a B-J equation with coefficients in ${\hat R}$.

\par \vskip 1pc \noindent
{\bf Step 1.} Solve $d_4 = 0$. We get the following expression of $u$,
which will be substituted to all of $d_i$:
$$u = \frac{1}{5}(-3p \sigma_3 + 2w \sigma_2 +
4q \sigma_4 - 2q \sigma_2^2).$$

\par \vskip 1pc \noindent
{\bf Step 2.} Note that $d_3$ is a quadratic form in $v, w, p, q$.
We find the standard normal form of $d_3$:
$$d_3 = \mu_1 \lambda_1^2 - \mu_2 \lambda_2^2 + \mu_3 \lambda_3^2 - \mu_4 \lambda_4^2,$$
where $\mu_i$ are elements in $R$ and
$\lambda_j$ are linear forms in $v, w, p, q$ with coefficients in $R$.
Write $\mu_2/\mu_1 = \mu_{2,1} (\mu_{2,1}')^2$ such that
$\mu_{2,1}$ is square free in $R$ and $\mu_{2,1}'$ is in the fraction
field of $R$. Then the extension ${\hat R}_1 := R[\sqrt{\mu_{2,1}}]$
of $R$ is a normal ring, since $R$ is \text{\rm UFD}
(see [Ma, \S 9, Example 4, p. 65]). So the singular locus
of ${\Cal Spec} \, {\hat R}_1$ is of co-dimension at least two, outside of which
${\hat R}_1$ is regular and hence a \text{\rm UFD}.

\par
Over the smooth locus of ${\Cal Spec} \, {\hat R}_1$, we write
$\mu_4/\mu_3 = \mu_{4,3}(\mu_{4,3}')^2$ where
$\mu_{4,3}$ is a square free regular function of ${\Cal Spec} \, {\hat R}_1$
and $\mu_{4,3}'$ is in the fraction
field of ${\hat R}_1$. Then the extension
${\hat R}_2 := {\hat R}_1[\sqrt{\mu_{4,3}}]$
of ${\hat R}_1$ is a normal ring in the open set of ${\Cal Spec} \, {\hat R}_2$
lying over the smooth locus of ${\Cal Spec} \, {\hat R}_1$. Hence
the singular locus of ${\Cal Spec} \, {\hat R}_2$ is of co-dimension at least two.

\par
Now any solution $(v, w, p, q)$ satisfying linear equations in $v, w, p, q$ below
will also satisfies the equation $d_3 = 0$
$$\lambda_1 - \sqrt{\mu_{2,1}} \mu_{2,1}'\lambda_2 = 0, \,\,\,\,
\lambda_3 - \sqrt{\mu_{4,3}} \mu_{4,3}'\lambda_4 = 0.$$

\par \vskip 1pc \noindent
{\bf Step 3.} Note that $d_2$ is a cubic form in $v, w, p, q$.
We substitute the two linear equations in Step 2 into $d_2$.
Then $d_2$ will become a cubic form in only two of the 4 variables $v, w, p, q$,
say in $w, q$ only. Consider the general case that
the Galois group of the cubic polynomial $d_2/q^3$ in $w/q$
over the fraction field of ${\hat R}_2$ is equal to $S_3$.
Taking a linear transformation
of coordinates $(w, q)$ over ${\hat R}_2$, we may
assume that $d_2$ multiplied by some non-zero elements in ${\hat R}_2$,
is equal to the following cubic form in $(w, q)$ over ${\hat R}_2$
(where the new $w$ and $q$ are ${\hat R}_2$-linear combination
of the old $w$ and $q$)
$$w^3 + s_1 w q^2 + t_1 q^3.$$
As in \S 1, if $s_1 \ne 0$, over the smooth locus of ${\Cal Spec} \, {\hat R}_2$, we
write $s_1 = a_{10}a_{11}a_{12}^2b_{11}, t_1 = b_{10}a_{11}a_{12}^2b_{11}^2$
and the discriminant $4s_1^3 + 27 t_1^2 = (a_{11}a_{12}^2)^2b_{11}^3 c_{11}c_{10}^2$. Let ${\hat R}_3 := {\hat R}_2[\sqrt{b_{11}c_{11}}]$, which is
normal over the smooth locus of ${\Cal Spec} \, {\hat R}_2$, so the singular locus
of ${\Cal Spec} \,{\hat R}_3$ is of co-dimension at least 2.

\par
Suppose that $(w, q)$ is a zero of the cubic form above. We now
define $\gamma$ (set $\gamma = w/q$ if $s_1 = 0$):
$$
\gamma := 6a_{10}(w/q)^2-9b_{10}b_{11}(w/q)+
\sqrt{3}c_{10}\sqrt{b_{11}c_{11}}(w/q)+4a_{10}^2a_{11}a_{12}^2b_{11}.
$$
Then, using the fact that
$c_{11}c_{10}^2 = 4a_{11}a_{12}^2a_{10}^3 + 27 b_{11}b_{10}^2$,
one can check that $\gamma$ satisfies an equation
$$\gamma^3 = \ell_1\ell_2^2\ell_3^3,$$
where $\text{\rm div} \,(\ell_i)$ ($i=1, 2$) are reduced divisors of ${\Cal Spec} \, {\hat R}_3$
and $\ell_3$ is in ${\hat R}_3$.
Since ${\hat R}_3[\gamma] \subset {\hat R}_3[w/q]$ and since these two rings have
the same degree over ${\hat R}_3$, they have the same fraction field
(= the splitting field of $d_2/q$ in $w/q$ over ${\hat R}_2$)
and the same normalization. We define
$$\align
{\hat R} &= {\hat R}_3 + {\hat R}_3 \gamma_1 + {\hat R}_3 \gamma_2, \\
\gamma_k &= \dfrac{\gamma^k}{\prod_{j=1}^3 \ell_j^{[jk/3]}}.
\endalign $$
By Lemma 1.3, ${\hat R}$ is a rank-3 free ${\hat R}_3$-module
and coincides with the normalization
of the ring ${\hat R}_3[\gamma]$ (or equivalently of the ring ${\hat R}[w/q]$)
on its open set lying over the smooth locus
of ${\Cal Spec} \, {\hat R}_3$.
So the singular locus of ${\Cal Spec} \, {\hat R}$ is of co-dimension at least 2.
Also for some $u_k$ in ${\hat R}_3$, we have the following,
since $w/q$ is integral over ${\hat R}_2$ (and also over ${\hat R}_3$)
$$\dfrac{w}{q} = \sum_{k=0}^2 u_k \gamma_k.$$

\par \vskip 1pc \noindent
{\bf Step 4.}
Using the above linear relation, the two linear equations at the end of Step 2
and the linear equation in Step 1, we can write
$$u = u_1q, \,\, v = v_1q, \,\,\, w = w_1q, \,\,\, p= p_1q$$
such that these 4 coefficients of $q$ are in the fraction field $Q({\hat R})$ of
${\hat R}$ and that
$y := u+v\alpha+w\alpha^2+p\alpha^3+q\alpha^4$
satisfies a type B-J equation : ${\hat y}^5 + \hat s {\hat y} + \hat t = 0$,
where $\hat s = d_1 = {\hat s}_1q^4, \hat t = d_0 = {\hat t}_1 q^5$
with coefficients of $q^4, q^5$ in $Q({\hat R})$.
Replacing $y$, ${\hat s}$, ${\hat t}$ by their multiples
of elements in ${\hat R}$, we may assume that $y$ is
in ${\hat R}[\alpha]$ already and satisfies a B-J equation defined over ${\hat R}$.

\par \vskip 1pc \noindent
{\bf (4.2).} Here are some detailed calculations which are followed by an example.
In Step 2 above, if $\sigma_2 \ne 0$ and if $\tau_i \ne 0$ ($i = 1, 2$), where
$$\align
\tau_1 = &45 \sigma_3^2+12\sigma_2^3-40\sigma_2 \sigma_4, \\
\tau_2 = &160\sigma_4^3+117\sigma_4\sigma_3^2\sigma_2+
12\sigma_2^4\sigma_4-88\sigma_2^2\sigma_4^2-4\sigma_2^3\sigma_3^2 + \\
         &(-27)\sigma_3^4+
        125\sigma_2\sigma_5^2-40\sigma_2^2\sigma_3\sigma_5-300\sigma_4\sigma_3\sigma_5,
\endalign$$
then we can write $d_3 = V-W+P-Q$, where
$$\align
V  =&\sigma_2\left(v-\frac{1}{2 \sigma_2}\left(5q \sigma_5+2p
\sigma_2^2-4p \sigma_4+
3w \sigma_3-5q \sigma_3 \sigma_2\right)\right)^2, \\
W  =&\frac{\tau_1}{20 \sigma_2}\left(w-\frac{1}{\tau_1} W'\right)^2, \\
W' =& 60p\sigma_4 \sigma_3 +
      8\sigma_2^2p\sigma_3 -75q\sigma_3 \sigma_5 +
      45q\sigma_3^2\sigma_2 + \\
    & (-50)p\sigma_2 \sigma_5 +
      12\sigma_2^4q-44\sigma_2^2q\sigma_4, \\
P  =& \frac{\tau_2}{\tau_1}\left(p - \frac{q}{2 \tau_2} P'\right)^2, \\
P'= &195\sigma_3^2\sigma_2 \sigma_5 -375\sigma_3 \sigma_5^2+
36\sigma_2^4\sigma_5 -4\sigma_4 \sigma_3 \sigma_2^3+ \\
    &48\sigma_4^2\sigma_3 \sigma_2 +400\sigma_5 \sigma_4^2
-260\sigma_2^2\sigma_5 \sigma_4 -27\sigma_3^3\sigma_4, \\
Q = &\left(\frac{1}{20 \tau_1 \tau_2} Q' -4\sigma_3^2\sigma_2+
4\sigma_3\sigma_5-\frac{12}{5}\sigma_4\sigma_2^2+
\frac{2}{5}\sigma_4^2+\frac{3}{5}\sigma_2^4\right)q^2.
\endalign$$
Here $Q'$ is a homogeneous polynomial of degree $26$ over ${\bold Z}$
in $\sigma_i$, where we set $\text{\rm deg}(\sigma_i) = i$.
In notation of Step 2,
$$\align
\dfrac{\mu_2}{\mu_1} &= \frac{\tau_1}{20 \sigma_2^2}, \\
\dfrac{\mu_4}{\mu_3} &= \frac{1}{20 \tau_2^2} \left(Q' + 20
\tau_1\tau_2\left(-4\sigma_3^2\sigma_2+
4\sigma_3\sigma_5-\frac{12}{5}\sigma_4\sigma_2^2+
\frac{2}{5}\sigma_4^2+\frac{3}{5}\sigma_2^4\right)\right).
\endalign$$

\par \vskip 1pc \noindent
{\bf Example 4.3.} We choose $\sigma_i$ below so that $\mu_{i+1}/\mu_i$
in {\bf (4.2)} above are relatively simpler:
$$\sigma_2 = \frac{3}{10}, \,\,\, \sigma_3 = \frac{1}{150}, \,\,\,
\sigma_4 = \frac{21}{2000}, \,\,\, \sigma_5 = \frac{-427}{75000}.$$
Then the linear equations mentioned in Step 4 above are:
$$\align
u &= \frac{1}{5}\left(-3 p \sigma_3 + 2 w \sigma_2 +
     4 q \sigma_4 - 2 q \sigma_2^2\right), \\
v &= \frac{1}{9000}\left(654p+3300w-1463q\right), \\
p &= \frac{1}{6}q, \\
L &= 0.
\endalign $$
Here $L$ is a linear factor of the cubic form below
(which is $d_2$ multiplied by $2^4 \times 3^7 \times 5^{11}$)
$$
612630271q^3+900w\left(-2004300qw-2283643q^2+7590000w^2\right).
$$

\par \vskip 1pc \noindent
{\bf (4.4).}
Now we shall calculate the integral closure of $R[\alpha]$.
Consider the general case where the polynomial
${\hat y}^5 + {\hat s} {\hat y} + {\hat t} \in \hat R[\alpha]$
found in {\bf (4.1)}, is irreducible over the fraction field of ${\hat R}$.
Let $\hat y = y \in R[\alpha]$ be a zero of this polynomial
as found in {\bf (4.1)}.

\par
Set ${\hat A} := {\hat R}[y]$.
Then ${\hat A}$ is a degree-5 extension of ${\hat R}$.
Since the extension ${\hat R}[\alpha]$ of ${\hat R}$ contains ${\hat A}$
and $z = \alpha$ satisfies the degree-5 polynomial $f(z)$
in $R[z] \subset {\hat R}[z]$, these two extensions of ${\hat R}$
have the same fraction field (and $f(z)$ is irreducible over
the fraction field of ${\hat R}$)
and hence the same normalization, which we denote by $A^*$.

\par
We denote by $A^*_0$ the ${\hat R}$-submodule of ${A^*}$
consisting of elements of trace zero over ${\hat R}$.
We will factor any element in $\hat R$ into the product of primes over the
smooth locus of ${\Cal Spec} \, \hat R$. We do this for $\hat s$ and $\hat t$,
the coefficients of the B-J polynomial in {\bf (4.1)}.
Note that $\hat a_i, \hat b_j$ below may have
different expressions at different affine open sets of ${\Cal Spec} \, \hat R$,
and $\text{\rm div}(\hat a_i)$ and $\text{\rm div}(\hat b_i)$ are reduced
divisors for all $i \ge 1$
$$
\hat s = \hat a_0 \, \hat a_1\,\hat a_2^2\,\hat a_3^3\,\hat a_4^4 \,
\hat b_1^3\,\hat b_2^2\, \hat b_3, \hskip0.3cm
\hat t = \hat b_0 \, \hat a_1\,\hat a_2^2\,\hat a_3^3\,\hat a_4^4 \,
\hat b_1^4\,\hat b_2^3\,\hat b_3^2.
$$
We can also define $\hat c_k, \hat f_k, \hat g_k, \hat h_k, \beta_k, \cdots$
over the smooth locus of ${\Cal Spec} \, \hat R$ as in \S 2
(with $R[\alpha]$ there replaced by ${\hat R}[y]$ here).
By Theorem 2.1, $A_0^*$ is given as follows
(over the smooth locus of ${\Cal Spec} \, \hat R$):
$$
A_0^* =\left\{\left.
\dfrac{\hat u_1\beta_1+\hat u_2
\beta_2+\hat u_3\beta_3 + \hat u_4 \beta_4}{\hat c_0}
\,\right|\, \hat u_j\in {\hat R}, \ \hat c_0\,|\,\hat f_k\hat u_k+\hat g_k\hat u_{k+1},
\ k=1, 2, 3 \,\right\}.
$$
In the expression above, ${\hat c}_0 | (\hat f_k\hat u_k+\hat g_k\hat u_{k+1})$
means that
${\hat f}_k {\hat u}_k + {\hat g}_k {\hat u}_{k+1} - {\hat c}_0 {\hat u}_{k+4} = 0$
for some ${\hat u}_{k+4}$ in ${\hat R}$.
Using bases, we have the expressions below, where $u_k$ are in ${\hat R}_3$,
$v_k$ are in ${\hat R}_2$, $w_k$ are in ${\hat R}_1$
and $x_k$ are in $R$
$$\align
{\hat u}_k &= u_k + u_{k+7} \gamma_1 + u_{k+14} \gamma_2 \,\,\, (1 \le k \le 7), \\
u_k &= v_k + v_{k+21} \sqrt{b_{11}c_{11}} \,\,\, (1 \le k \le 21), \\
v_k &= w_k + w_{k+42} \sqrt{\mu_{4,3}} \,\,\, (1 \le k \le 42), \\
w_k &= x_k + x_{k+84} \sqrt{\mu_{2,1}}  \,\,\, (1 \le k \le 84).
\endalign $$
So we can write ${\hat u}_k = \sum_{j=1}^{168} \varepsilon_{kj} x_j$,
where $\varepsilon_{kj}$ is in the ring ${\bold Z}[\sqrt{\mu_{2,1}}, \sqrt{\mu_{4,3}},
\sqrt{b_{11}c_{11}}, \gamma_1, \gamma_2]$.
Now we can express a general element $x$ of $A_0^*$ as follows
$$\align
x &= \frac{\hat u_1\beta_1+\hat u_2
\beta_2+\hat u_3\beta_3 + \hat u_4 \beta_4}{\hat c_0} = \sum_{j=1}^{168} x_j d_j, \\
d_j &= \frac{\sum_{k=1}^4 \varepsilon_{kj} \beta_k }{\hat c_0}.
\endalign $$
We can calculate all $d_j$ ($1 \le j \le 168$) more explicitly
as follows, where those $d_j$ not listed below, are equal to 0
$$\align
{\hat c}_0 d_k &= \beta_k \,\,\, (1 \le k \le 4), \\
{\hat c}_0 d_k &= \gamma_1 \beta_{k-7} \,\,\, (8 \le k \le 11), \\
{\hat c}_0 d_k &= \gamma_2 \beta_{k-14}  \,\,\, (15 \le k \le 18),\\
{\hat c}_0 d_k &= \sqrt{b_{11}c_{11}}\beta_{k-21} \,\,\, (22 \le k \le 25), \\
{\hat c}_0 d_k &=\gamma_1 \sqrt{b_{11}c_{11}}\beta_{k-28} \,\,\, (29 \le k \le 32), \\
{\hat c}_0 d_k &= \gamma_2 \sqrt{b_{11}c_{11}} \beta_{k-35} \,\,\, (36 \le k \le 39),\\
{\hat c}_0 d_k &= \sqrt{\mu_{4,3}}\beta_{k-42} \,\,\, (43 \le k \le 46), \\
{\hat c}_0 d_k &= \gamma_1 \sqrt{\mu_{4,3}}\beta_{k-49} \,\,\, (50 \le k \le 53),\\
{\hat c}_0 d_k &= \gamma_2 \sqrt{\mu_{4,3}} \beta_{k-56} \,\,\, (57 \le k \le 60), \\
{\hat c}_0 d_k &= \sqrt{b_{11}c_{11}\mu_{4,3}}\beta_{k-63} \,\,\, (64 \le k \le 67), \\
{\hat c}_0 d_k &= \gamma_1 \sqrt{b_{11}c_{11}\mu_{4,3}}\beta_{k-70} \,\,\,
(71 \le k \le 74),\\
{\hat c}_0 d_k &= \gamma_2 \sqrt{b_{11}c_{11}\mu_{4,3}} \beta_{k-77} \,\,\,
(78 \le k \le 81), \\
{\hat c}_0 d_k &= \sqrt{\mu_{2,1}}\beta_{k-84} \,\,\, (85 \le k \le 88), \\
{\hat c}_0 d_k &= \gamma_1 \sqrt{\mu_{2,1}}\beta_{k-91} \,\,\,(92 \le k \le 95),\\
{\hat c}_0 d_k &= \gamma_2 \sqrt{\mu_{2,1}} \beta_{k-98} \,\,\, (99 \le k \le 102), \\
{\hat c}_0 d_k &= \sqrt{b_{11}c_{11}\mu_{2,1}}\beta_{k-105}\,\,\,(106 \le k \le 109),\\
{\hat c}_0 d_k &= \gamma_1 \sqrt{b_{11}c_{11}\mu_{2,1}}\beta_{k-112} \,\,\,
(113 \le k \le 116),\\
{\hat c}_0 d_k &= \gamma_2 \sqrt{b_{11}c_{11}\mu_{2,1}} \beta_{k-119} \,\,\,
(120 \le k \le 123), \\
{\hat c}_0 d_k &= \sqrt{\mu_{4,3}\mu_{2,1}}\beta_{k-126}\,\,\,(127 \le k \le 130),\\
{\hat c}_0 d_k &= \gamma_1 \sqrt{\mu_{4,3}\mu_{2,1}}\beta_{k-133} \,\,\,
(134 \le k \le 137),\\
{\hat c}_0 d_k &= \gamma_2 \sqrt{\mu_{4,3}\mu_{2,1}} \beta_{k-140} \,\,\,
(141 \le k \le 144), \\
{\hat c}_0 d_k &= \sqrt{b_{11}c_{11}\mu_{4,3}\mu_{2,1}}\beta_{k-147} \,\,\,(148 \le k \le 151),\\
{\hat c}_0 d_k &= \gamma_1 \sqrt{b_{11}c_{11}\mu_{4,3}\mu_{2,1}}\beta_{k-154} \,\,\,
(155 \le k \le 158), \\
{\hat c}_0 d_k &= \gamma_2 \sqrt{b_{11}c_{11}\mu_{4,3}\mu_{2,1}} \beta_{k-161} \,\,\,
(162 \le k \le 165).
\endalign$$

\par
Now we look at the syzygy condition in $A_0^*$.
We can write
$$\gather
{\hat f}_k {\hat u}_k + {\hat g}_k {\hat u}_{k+1} -
{\hat c}_0 {\hat u}_{k+4} = \sum_{j=1}^{168} x_j e_{kj}, \\
e_{kj} = \hat f_k \varepsilon_{kj} + \hat g_k \varepsilon_{k+1,j} -
\hat c_0 \varepsilon_{k+4, j}.
\endgather $$

\par
Let $A := R[\alpha]$, let ${\wt A}$ be the normalization of
$A$ in the fraction field of $A$ and let
${\wt A}_0$ be the $R$-submodule of ${\wt A}$
consisting of trace (over $R$) free elements.
Clearly, ${\wt A}_0 = \{\text{\rm Tr}(x) \,| \,
x \in A_0^*\}$ (see {\bf (4.5)} and {\bf (4.6) (2)} below for Tr).
Thus we have proved the main theorem below of this section,
since an element in a normal ring ${\wt A}$
is determined uniquely by its restriction to the complement of
a co-dimension 2 subset.

\proclaim{Theorem 4.5}
Let $R$ be a Noetherian \text{\rm UFD} containing a field with $\text{\rm Char} \, R$ coprime to $2, 3, 5$.
For a general degree-$5$ extension $R \subset R[\alpha]$,
where $\alpha$ is a root of an irreducible polynomial
$f(z) = z^5+\sigma_2 z^3 -\sigma_3 z^2+\sigma_4 z - \sigma_5$
in $R[z]$, the integral closure $\widetilde A$ of $R[\alpha]$ in the fraction
field of $R[\alpha]$ is given by ${\wt A} = R \oplus {\wt A}_0$, and
the trace (over $R$) free $R$-submodule $\widetilde A_0$ of
${\wt A}$ is given as follows:
$${\wt A}_0 = \left\{\sum_{i=1}^{168} x_i \text{\rm Tr}(d_i) \ |
\ x_i \in R, \ \sum_{j=1}^{168} x_j e_{kj} = 0, \ \
(k=1,2,3)\right\},$$ where $\text{\rm Tr}$ is the trace for the
field extension $Q(R[\alpha]) \subset Q({\hat R}[\alpha])$.
\endproclaim

\remark{Remark 4.6} $(1)$ In Theorem $4.5$, by a general
degree-$5$ extension, we mean that the conditions that $\sigma_2
\ne 0$ and $\tau_i \ne 0$ ($i = 1, 2$) in {\bf (4.2)} are
satisfied and that the extension $R \subset {\hat R}$ is of degree
$24$. If the extension has a smaller degree, the same process of
reducing to a B-J extension works and is even simpler by a least
one step.

$(2)$ For the calculation of $\text{\rm Tr}(d_i)$, we note that
$$\text{\rm Tr}|\hat R[\alpha]/R[\alpha]
= \hat T_1 \circ \hat T_2 \circ \hat T_3 \circ \hat T_4,$$
and $\hat T_i = \text{\rm Tr} | \hat R_i[\alpha]/\hat R_{i-1}[\alpha]$
is the lifting of
$T_i = \text{\rm Tr} | \hat R_i / \hat R_{i-1}$,
where $\hat R_0 := R$ and $\hat R_4 := \hat R$.
The traces $T_i$ ($i = 1, 2, 3$) for quadratic extensions are rather easy
and for $T_4$ we have $T_4(u_0 + u_1 \gamma_1 + u_2 \gamma_2) = u_0$
whenever $u_i \in \hat R_3$
\endremark

\head
5. Some applications in algebraic geometry
\endhead

We will now apply the calculation of integral closure
to B-J covers of a factorial variety.

\par
Let $X$ be a factorial variety over a field
$k$ with $\text{\rm Char} \, k$ coprime to both $n-1$ and $n$, let $L$
be a line bundle over $X$, and let $s$ and $t$ be two non-zero
global sections of $L^{n-1}$ and $L^n$, respectively.
Then we can construct a normal finite cover $\pi: Y\to X$ of
degree $n$ by adding a root of the B-J polynomial $f=z^n+sz+t$,
which is irreducible over the function field $k(X)$.
We may also assume that the data $(s, t)$ is minimal.
The construction is as follows.

\par
Let $p: [L]\to X$ be the $\bold A^1$-bundle over $X$ associated with the
line bundle $L$. Denote by $z$ the fiber coordinate. Then
$f=z^n+sz+t$ is a global section of $p^*(L^n)$. Denote by $\overline{Y}$
the zero scheme of $f$. Let $Y$ be the normalization of $\overline{Y}$.
Then we see that the induced cover $\pi:Y \to X$ is a finite
morphism of degree $n$. We call it a B-J {\it cover}.

\par
We have the same type of factorizations of $s$ and $t$ as in \S1.
So we can define global sections $f_i$, $g_i$ and $h_i$ ....
of some line bundles as in \S 2. The globalization of Proposition 1.2 shows that
there is a one to one correspondence between the B-J polynomials
$z^n+sz+t$ with the data $(s, t)$ minimal
and the triplets $(a,b,c)$ of coprime sections of a line
bundle with $a+b=c$. So $a+b=c$ can be viewed as the covering data
of a B-J cover.

\par
We denote by $E_\pi$ the trace free subsheaf of $\pi_*\Cal  O_Y$.
Since $\text{\rm Char} \,  k$ is coprime to $n$, the
trace map $\text{Tr}: \pi_*\Cal  O_Y \to \Cal  O_X$ splits
$\pi_*\Cal  O_X$ as the direct sum of $\Cal O_X$ and $E_\pi$:
$$
\pi_*\Cal  O_Y =\Cal O_X\oplus E_\pi.
$$

\par
Now we construct a $(n-2)\times (2n-3)$ matrix $M=(m_{ij})$:
$$
m_{ij}:=\cases f_i, &\text{ if } \ j=i, \cr
g_i, &\text{ if } \ j=i+1, \cr
c_0, &\text{ if } \ j=i+n-1,\cr
0, &\text{ otherwise}.
\endcases
$$
Define $2n-3$ hypersurfaces ($g_0:=1, \ f_{n-1}:=1$):
$$
V_i:=\cases
\text{div}(g_1\cdots g_{i-1}f_i\cdots f_{n-2}),
&\text{ if } \ 1\leq i\leq n-1,\cr
\text{div}(c_0g_1\cdots g_{k}f_{k+2}\cdots f_{n-2}),
&\text{ if } \ i=n+k, \ 0\leq k \leq n-3.\cr
\endcases
$$
$V_1,\cdots,V_{2n-3}$ determine a syzygy sheaf $\Cal F$
$$\align
0\to \Cal F\to &\oplus_{i=1}^{2n-3}\Cal O_X(-V_i) \ {\overset M
\to\longrightarrow} \ \oplus_{i=1}^{n-2}\Cal O_X(-V_i+\text{div}(f_i)),\cr
&(v_1,\cdots,v_{2n-3})^t\mapsto M(v_1,\cdots,v_{2n-3})^t.
\endalign
$$
We have chosen $V_i$ so that the map $M$ is well defined. Namely,
the $i$-th component of an element in the image of $M$ is a ${\bold Z}$-combination
of local sections in the same line bundle ${\Cal O}_X(-V_i + \text{\rm div}(f_i))$
with a fixed transition function.
One can check easily that the divisor
$$
T:=V_i-iL+\text{div}(c_0)+\text{div}(h_i)
$$
is independent of $i = 1, \dots, n-1$.
Now we state the main result of the section:

\proclaim{Theorem 5.1} Let $\pi : Y \rightarrow X$ be a degree-$n$
finite morphism from a normal variety onto the factorial variety
$X$ defined over a field $k$ with $\text{\rm Char} \, k$ coprime
to $n$ and $n-1$, so that the function field $k(Y)$ is obtained by
adding to $k(X)$ a root of a polynomial $f = z^n + s z + t$, where
$0 \ne s \in H^0(X, L^{n-1})$, $t \in H^0(X, L^n)$ for a line
bundle $L$ with the data $(s, t)$ minimal in the sense of ${\bold
(1.1)}$ $($for $\pi$, see the beginning of \S 5 for more
detail$)$. Then $\pi_*\Cal  O_Y =\Cal O_X\oplus E_\pi$ and the
trace $($over $X)$ free part $E_{\pi} \cong \Cal F(T)$.
\endproclaim

\demo{Proof} We define a map
$$
\tau: \Cal F(T)\to E_\pi, \hskip0.3cm
(v_1,\cdots,v_{2n-3})^t\mapsto {\sum_{i=1}^{n-1}v_i\beta_i\over c_0}.
$$
Since $\beta_i/c_0=\alpha^i/h_ic_0$ can be viewed as a basis of $\Cal O_X(T-V_i)$,
we can check that the linear map is well defined globally.
By Theorem 2.1, $\tau$ is an isomorphism locally and hence globally.
\enddemo

\par
Note that even if $X$ is only normal, the above argument works
too. In fact, the above study works for the induced B-J cover
over the smooth part of $X$. Since
the singular locus of $X$ is of co-dimension at least 2, most
of the results can be extended to $X$.

\par
If $X$ is a smooth surface, then the singularities of any B-J
cover $\pi:Y\to X$ can be resolved by the following classical
method. Let $a+b=c$ be the data of $\pi$. If the branch locus has
a singular point $p_1$ worse than normal crossing, then we blow
up $X$ at $p_1$: $\sigma_1:X_1\to X$. Let $\pi_1:Y_1\to X_1$
be the pullback of $\pi$ via $\sigma_1$, i.e., $Y_1$ is the
normalization of $X_1\times _XY$. This $\pi_1$ is a B-J
cover defined by $z^n+\sigma_1^*(s)z+\sigma_1^*(t)=0$.
If we denote by
$a^{(1)}+b^{(1)} = c^{(1)}$ the minimal data of $\pi_1$, then
it is obtained from
$$
\sigma_1^*(a)+\sigma_1^*(b)=\sigma_1^*(c)
$$
by eliminating the common factors from both sides, which come
from the exceptional divisor. After a finite number of steps,
the branch locus of $\pi_m$ becomes normal crossing. Now we can
resolve the singularities of $Y_m$ by the Hirzebruch-Jung
method.

\par
When $n=3$, it has been proved in [Ta2] that the
singularities of $Y_m$ can be resolved by simply applying
the same process as above to the singular points of the branch locus.
Thus for any triple cover, we have a ``canonical resolution''
of the singularities the same as in the double cover case (see [Ta2, Theorem 7.2]).

\par \vskip 1pc
The result below follows from Theorem 5.1, the fact that $(n-1)L
\sim \text{\rm div}(s)$ and the fact that the cokernel of $M$ in
the exact sequence preceding Theorem 5.1 is a direct sum of
sheaves each summand of which is supported on a set $\{f_i = g_i =
c_0 = 0\}$ (which has codimension at least 2); see Cor. 2.2 for
the local expression of $M$.

\proclaim{Corollary 5.2}
{\it With the same assumptions as in Theorem $5.1$,
we have linear equivalences:}
$$c_1(\pi_*{\Cal O}_Y) \sim (n-1)T - \sum_{i=n-1}^{2n-3} V_i - \sum_{i=1}^{n-2}
\text{\rm div}(f_i) \sim \text{\rm div}
\left(\frac{c_0}{s}\prod_{i=1}^{n-2}
\left(\frac{f_i}{g_i}\right)^{n-1-i}\right).$$
\endproclaim

\par
For a divisor $D$ on a variety, we denote by $D_{\text{\rm red}}$
the reduced divisor with the same support as $D$.
For a morphism $f$ we denote by $R_f$ the
ramification divisor and $B_f = f_*(R_f)$ the branch locus of $f$;
when $f$ is a finite morphism of degree $n$, we
write $R_f = R_{fs} + R_{ft} + \cdots$, $B_{fs} := f_*R_{fs}$, $B_{ft} := f_*R_{ft}$,
where $f$ is {\it simply} (resp. {\it totally}) ramified along $R_{fs}$
(resp. $R_{ft}$), i.e., over a generic point of $(B_{fs})_{\text{\rm red}}$
(resp. $(B_{ft})_{\text{\rm red}}$) the ramification index(es)
are $2, 1, \dots, 1$ (resp. is $n$). So
$R_{ft} = (n-1)(R_{ft})_{\text{\rm red}}$,
$B_{ft} = (n-1)(B_{ft})_{\text{\rm red}}$;
$R_{fs}$ and $B_{fs}$ are reduced.
When $n = 3$, we have $R_f = R_{fs} + R_{ft}$.

\proclaim{Theorem 5.3}  Let $\pi : Y \rightarrow X$ be a
degree-$3$ finite morphism from a normal variety onto the
factorial variety $X$ defined over a field $k$ with $\text{\rm
Char} \, k \ne 2, 3$. Then we have:
 \roster
 \item  We have a
linear equivalence:
$$c_1(\pi_*{\Cal O}_Y) \sim -(B_{\pi t})_{\text{\rm red}} - \eta,$$
where $\eta$ is a divisor satisfying $2\eta \sim B_{\pi s}$.

\item $$2 c_1(\pi_*{\Cal O}_Y) \sim - B_{\pi}.$$

\item {\it $\pi$ is unramified outside a co-dimension two
subset if and only if $2 c_1(\pi_*{\Cal O}_Y) \sim 0$.}
\endroster\endproclaim

\par \vskip 1pc
We now prove Theorem 5.3.
Since $\text{\rm Char} \, k \ne 3$, we see that $\pi$
is a B-J cover and given by a polynomial below (see [Ta2, Theorem 7.2]):
$$f(z) = z^3 + s z + t$$
where $s, t$ are global sections of $L^2$, $L^3$, where
$L$ is a line bundle on $X$. Since $f(z)$ is defined globally,
$z$ (a zero of $f$), $s$, $t$ have the transition functions
$\sigma_{ij}, \sigma_{ij}^2, \sigma_{ij}^3$, respectively, with respect to affine open sets
$\{U_i\}$ covering $X$.
We may assume that the data $(s, t)$ is minimal.
Thus Theorem 5.3 follows from Theorem 1.4 and Lemma 5.4 below.

\par
We will use the notations $A_i = \text{\rm div}(a_i), B_i = \text{\rm div}(b_i),
C_i = \text{\rm div}(c_i)$; see \S 1 for $a_k$, $b_k$ and $c_k$.
Lemma 5.4, where (1) was proved in [Ta1, Ta2], is now a consequence of Corollary 5.2;
in applying, we also use the fact that
$\text{\rm div}(s^{n}/t^{n-1}) \sim 0$ and
$C_1 + 2C_0 \sim B_1 + 2B_0$ from the definition of $c_i$ (for {\bf (5.4) (1)}).

\proclaim{Lemma 5.4} We assume the hypothesis and notation in Theorem $5.1$ or Corollary $5.2$.
 \roster
 \item {\it
Let $n = 3$, i.e., let $\pi : Y \rightarrow X$ be a degree $3$
finite morphism from a normal variety onto the factorial variety
$X$ defined over a field $k$ with $\text{\rm Char} \, k \ne 2, 3$,
so that the function field $k(Y)$ is obtained by adding to $k(X)$
a root of a polynomial $f = z^3 + s z + t$, where $0 \ne s \in
H^0(X, L^{2})$, $t \in H^0(X, L^3)$ for a line bundle $L$ with the
data $(s, t)$ minimal in the sense of ${\bold (1.1)}$ {$($}for
$\pi$, see the beginning of \S 5 for more detail$)$. Then we have:
$$c_1(\pi_*{\Cal O}_Y) = C_0 - A_1 - A_2 - B_0 - B_1 = -(A_1+A_2) - \eta,$$
where $\eta = B_0+B_1-C_0$ and $2\eta \sim B_1 + C_1$.}

\item {\it If $n = 4$, then we have the
linear equivalence:
$$c_1(\pi_*{\Cal O}_Y) \sim \text{\rm div}\left(\frac{c_0}{a_0^2a_1^2a_2^2a_3^3b_1b_2}\right).$$

\item {\it If $n = 5$, then we have:
$$c_1(\pi_*{\Cal O}_Y) \sim \text{\rm div}\left(\frac{c_0}{\left(a_1a_2a_3a_4b_0b_1b_2\right)^2b_3^3}\right).$$
\endroster
\endproclaim

\par \vskip 1pc
The first attempted proof of the following result appeared in [Ta2].
The approach here is different.

\proclaim{Theorem 5.5}
{\it Let $\pi : Y \rightarrow X$ be a finite morphism
of degree $3$ from a normal variety onto a factorial variety
defined over a field $k$ with
$\text{\rm Char} \,  k \ne 2, 3$. Suppose that $H^0(X, {\Cal O}_X) = k$
and every element of $k$ has a square root in $k$.
Then $\pi$ is Galois if and only if the
following two conditions are satisfied.}

\roster\item {\it Outside a co-dimension $2$ subset, $\pi$ is
either unramified or totally ramified, i.e., either $B_{\pi} = 0$
or $\pi^*((B_{\pi})_{\text{\rm red}}) = 3(R_{\pi})_{\text{\rm
red}}$.}

\item {\it $c_1(\pi_*{\Cal O}_Y) \sim - (B_{\pi})_{\text{\rm
red}}$.} \endroster
\endproclaim

\demo{Proof} We may assume that $\pi$ is given as in Lemma 5.4
[Ta2, Theorem 7.2], though $s$ might be identically 0.

\par
Suppose that $\pi$ is Galois. Then (1) is clear and we may assume that
$f = z^3 + \ell_1\ell_2^2$ with $\text{\rm div}(\ell_i)$ reduced.
By Lemma 1.3, we have $\pi_*{\Cal O}_Y =
{\Cal O}_X \oplus {\Cal O}_X z \oplus {\Cal O}_X (z^2/\ell_2)$.
Hence we have:
$$c_1(\pi_*{\Cal O}_Y) = -\text{\rm div}(z) - \text{\rm div}(z^2/\ell_2)
= - \text{\rm div}(z^3/\ell_2) = - \text{\rm div}(\ell_1\ell_2^2/\ell_2)
= - (B_{\pi})_{\text{\rm red}}.$$

\par
Now we assume (1) and (2). In notation of Lemma 5.4, $\pi$ is Galois
if $s = 0$. So we may assume that $s \ne 0$ and all conditions of Lemma 5.4
are satisfied.
As in \S 1, we have the expression below,
where $a_k, b_k, c_k$ are global sections of some line bundles
$$s = a_0a_1a_2^2b_1, \,\,\,\,\,\, t = b_0a_1a_2^2b_1^2.$$
Note that the discriminant $\delta$ of $f(z)$ is given by
$$\delta = (a_1a_2^2)^2b_1^3c_1c_0^2, \,\,\,
c_1c_0^2 = 4a_1a_2^2a_0^3 + 27b_1b_0^2.$$
Write $a_k = \{a_k(i)\}, b_k = \{b_k(i)\}, c_k = \{c_k(i)\}$
where $a_k(i), ...$ are regular sections on the affine open set $U_i$,
the union of which covers $X$.
From the abstract algebra, we know that when $\text{\rm Char} \, k \ne 2$,
the finite morphism $\pi : Y \rightarrow X$ is Galois
if and only if $\delta$ has a square root in the function field $k(X)$, i.e.,
$b_1(i)c_1(i)$ has a square root in $k(X)$ for one $i$ (and hence for all $i$).

\par
Note that $B_1 = C_1 = 0$
for otherwise $\pi$ would be of simple ramification over
$B_1 = \text{\rm div}(b_1)$ and $C_1 = \text{\rm div}(c_1)$ (Theorem 1.4),
contradicting (1).
So $b_1(i)$ and $c_1(i)$ are invertible on $U_i$.
We may assume that $b_1(i) = 1$ by renaming $a_0b_1$ and $b_0b_1^2$
as new $a_0$ and $b_0$, respectively so that $s = a_0a_1a_2^2$ and $t = b_0a_1a_2$.
Note that for all $i$, we have
$$c_1(i)c_0(i)^2 = 4a_1(i)a_2(i)^2a_0(i)^3 + 27b_0(i)^2.$$
The 3 terms in the equation
have the same transition function; indeed this equation is deduced by
renaming the quotient of $\delta = 4s^3 + 27t^2$ after the division by
the common factors of the two terms in $\delta$, while
the transition functions of $s, t$ are $\sigma_{ij}^2, \sigma_{ij}^3$, respectively
(so the two terms in $\delta$ have the same transition function).

\par
On the other hand, from the condition (2) and Lemma 5.4, we see that
$B_0 \sim C_0$. We may assume that $b_0 = \{b_0(i)\}$ and $c_0 = \{c_0(i)\}$
have the same transition function after adjusting $c_1$ if necessary.
This and the similar fact after the display in the previous paragraph
force $c_1 = \{c_1(i)\}$ to have the constant $1$ as its transition function,
i.e., $c_1(i) = c_1(j)$ for all $i, j$. So $c_1(i)$ is a global invertible
function. Hence $c_1(i)$ is a constant in $k$ because $H^0(X, {\Cal O}_X) = k$.
Now $b_1(i)c_1(i) = c_1(i)$ has a square root in $k \subset k(X)$.
So the Galoisness follows. This proves the theorem.
\enddemo

\remark{Remark 5.6} $(1)$ The condition $(1)$ alone in Theorem
$5.5$ is not enough to imply the Galoisness of $\pi$. Indeed, by a
Theorem of J. P. Serre, we know that for any $m \ge 2$ there is an
$m$-dimensional projective complex manifold $X$ with $\pi_1(X) =
S_3$, the symmetric group in 3 letters. 
Let $U$ be the universal cover of $X$ and let $Y = U/\langle
\sigma \rangle$ where $\sigma$ is an involution in $S_3$. Then $Y
\rightarrow X = U/S_3$ is a finite morphism of degree 3 between
smooth projective manifolds, which is unramified but non-Galois.

(2) Despite what we said in (1), the condition (1) in Theorem 5.5 together
with a condition (2)' that $Pic X$ has no 2-torsion element (this is true
when $\pi_1(X)$ has no index-2 subgroup) will imply the Galoisness of $\pi$.
Indeed, the conditions (1) and (2)' imply that $\eta = 0$ in Lemma 5.4
and hence the condition (2) in Theorem 5.5 holds.

$(3)$ We like to have a similar geometric Galoisness criterion for
degree 4, or 5 extension, but we realized from the discussion with
Professor Catanese that it is much more complicated. The following
is a partial result.

 $(4)$ In Theorem $5.5$ and Proposition $5.7
(2)$, reading the proof, we see that the two conditions that
$H^0(X, {\Cal O}_X) = k$ and every element in $k$ has a square
root in $k$ can be weakened to one condition that every element of
$H^0(X, {\Cal O}_X)$ has a square root in the function field
$k(X)$.
\endremark

\proclaim{\bf Proposition 5.7}
With the assumption in Theorem $5.1$, we set $n = 5$. Then we have:

$(1)$ If $\pi$ is Galois, i.e., if the Galois group $\text{\rm
Gal}(f)$ of the polynomial $f(z)$ is ${\bold Z}/(5)$, then the
following are true: \roster
\item"$(1a)$" $c_1(\pi_*{\Cal O}_Y) \sim
-2(B_{\pi})_{\text{\rm red}}$, and

\item"$(1b)$" Outside a co-dimension $2$ subset,
$\pi$ is either unramifield or totally ramified, i.e., either $B_{\pi} = 0$
or $\pi^*((B_{\pi})_{\text{\rm red}}) = 5(R_{\pi})_{\text{\rm red}}$.
\endroster

$(2)$ Conversely, suppose that $(1a)$ and $(1b)$ are satisfied and
suppose further that $H^0(X, {\Cal O}_X)$ $= k$ and every element
in the ground field $k$ has a square root in $k$. Then $\text{\rm
Gal}(f)$ is one of ${\bold Z}/(5)$, $D_{10}$ $($the dihedral group
of order $10)$ and $A_5$ $($the alternating group in $5$
letters$)$.
\endproclaim

\demo{Proof} Assume first that $\pi$ is Galois. Then $(1b)$ is
clear. We may also assume that the function field $k(Y)$ is
obtained by adding to $k(X)$ a root of the polynomial below
$$z^5 + \ell_1 \ell_2^2 \ell_2^3 \ell_4^4$$
where $\ell_j$ are reduced global sections of some line bundles.
Then, by Lemma 1.3, we have:
$$\pi_*{\Cal O}_Y = {\Cal O}_X \oplus {\Cal O}_X z \oplus
{\Cal O}_X (z^2/(\ell_3\ell_4)) \oplus {\Cal O}_X (z^3/(\ell_2\ell_3\ell_4^2))
\oplus {\Cal O}_X (z^4/(\ell_2\ell_3^2\ell_4^3)).$$
As in Theorem 5.5, we can now imply $(1a)$:
$$c_1(\pi_*{\Cal O}_Y) = \text{\rm div}(\ell_2^2\ell_3^4\ell_4^6/z^{10}) \sim -2 \text{\rm div}(\ell_1 \ell_2 \ell_3 \ell_4) =
-2(B_{\pi})_{\text{\rm red}}.$$
Here we used the fact that
$\text{\rm div}(z^{10}) \sim \text{\rm div}(\ell_1 \ell_2^2 \ell_2^3 \ell_4^4)^2$.

\par
Now assume $(1a)$ and $(1b)$. So $C_1 = \text{\rm div}(c_1) = 0$ and
$B_i = \text{\rm div}(b_i) = 0$ for all $i \ge 1$ (see Theorem
1.4). By Lemma 5.4, we have $c_1(\pi_*{\Cal O}_Y) \sim C_0 - 2(A_1
+ A_2 + A_3 + A_4 + B_0)$. This and $(1a)$ imply that $C_0 \sim
2B_0$. The rest is similar to Theorem 5.5. To be precise, since
$B_i = \text{\rm div}(b_i) = 0$ ($i \ge 1$) we can adjust $a_0$
and $b_0$ and assume that $b_i = 1$ ($i \ge 1$) so that $s =
a_0a_1a_2^2a_3^3a_4^4$ and $t = b_0a_1a_2^2a_3^3a_4^4$. So the
discriminant $\delta$ of the polynomial $f(z)$ in Theorem 5.1
equals $(a_1a_2^2a_3^3a_4^4)^4 c_1^2 c_0$ where
$$c_1c_0^2 = 256a_1a_2^2a_3^3a_4^4a_0^5+3125b_0^4.$$
The $3$ terms in the equation above have the same transition function;
since $C_0 \sim 2B_0$, we can adjust $c_1$ and assume that
$C_0 = \text{\rm div}(c_0)$ and $2B_0 = \text{\rm div}(b_0^2)$
have the same transition function. These two facts
imply that $c_1$ has constant 1 as its transition function, so $c_1$
is a non-zero scalar in $k$. Thus the discriminant $\delta$ of $f(z)$
has a square root in $k(X)$, whence $\text{\rm Gal}(f) \le A_5$
(noting that $\text{\rm Char} \, k \ne 2$ for it is coprime to $5-1$ (and also $5$)).
Now the assertion (2) follows from the classification of subgroups in $A_5$.
\enddemo

\enddocument